\documentclass[a4paper,12pt]{amsart}
\usepackage{ifthen}
\usepackage{mathrsfs}
\nonstopmode
\numberwithin{equation}{section}
\newtheorem{thm}{Theorem}[section]
\newtheorem{cor}[thm]{Corollary}
\newtheorem{lem}[thm]{Lemma}

\theoremstyle{definition}

\newenvironment{rem}{%
\bigskip
\noindent
\textsl{{\sl Remark. }}}{\bigskip}

\newenvironment{pf}[1][]{%
 \vskip 3mm
 \noindent
 \ifthenelse{\equal{#1}{}}%
  {{\slshape Proof. }}%
  {{\slshape #1.} }%
 }%
{\qed\bigskip}

\newcounter{alphabet}
\newcounter{tmp}
\newenvironment{Thm}[1][]{\refstepcounter{alphabet}%
\bigskip%
\noindent%
{\bf Theorem \Alph{alphabet}}%
\ifthenelse{\equal{#1}{}}{}{ (#1)}%
{\bf .}
\itshape}{\vskip 8pt}

\newcommand{\A}{{\mathcal A}}

\newcommand{\C}{{\mathbb C}}
\newcommand{\D}{{\mathbb D}}

\newcommand{\dom}{{\mathcal D}}

\newcommand{\K}{{\mathcal K}}
\newcommand{\M}{{\mathcal M}}

\newcommand{\N}{{\mathcal N}}

\newcommand{\rg}{{\mathcal R}}
\newcommand{\es}{{\mathcal S}}

\newcommand{\V}{{\mathcal V}}

\newcommand{\spl}{{\mathcal{SP}}}

\newcommand{\LU}{{\mathcal{LU}}}
\newcommand{\ZF}{{\mathcal{ZF}}}

\newcommand{\CC}{{\mathcal C}}
\newcommand{\st}{{\mathcal{SS}}}
\newcommand{\bD}{{\overline{\mathbb D}}}

\renewcommand{\Im}{\,{\operatorname{Im}\,}}

\renewcommand{\Re}{{\operatorname{Re}\,}}

\newcommand{\id}{{\operatorname{id}}}

\newcommand{\inv}{^{-1}}

\renewcommand{\arg}{\,{\operatorname{arg}\,}}

\newcommand{\Log}{{\operatorname{Log}\,}}
\newcommand{\Arg}{{\operatorname{Arg}\,}}

\newcounter{minutes}
\divide\time by 60
\newcounter{hours}
\multiply\time by 60
\addtocounter{minutes}{-\time}

\begin{document}
\bibliographystyle{amsplain}
\title{
On power deformations of univalent functions
}


\author[Y.~C.~Kim]{Yong Chan Kim}
\address{Department of Mathematics Education, Yeungnam University, 214-1 Daedong
Gyongsan 712-749, Korea}
\email{kimyc@ynu.ac.kr}
\author[T. Sugawa]{Toshiyuki Sugawa}
\address{Division of Mathematics, Graduate School of Information Sciences,
Tohoku University, Aoba-ku, Sendai 980-8579, Japan}
\email{sugawa@math.is.tohoku.ac.jp}
\keywords{univalent function, variability region, spirallike function}
\subjclass[2000]{Primary 30C45; Secondary 30C55}
\begin{abstract}
For an analytic function $f(z)$ on the unit disk $|z|<1$ with
$f(0)=f'(0)-1=0$ and $f(z)\ne0, 0<|z|<1,$ we consider the power deformation
$f_c(z)=z(f(z)/z)^c$ for a complex number $c.$
We determine those values $c$ for which the operator $f\mapsto f_c$
maps a specified class of univalent functions into the class of 
univalent functions.
A little surprisingly, we will see that 
the set is described by the variability
region of the quantity $zf'(z)/f(z),~|z|<1,$ for the class in most cases
which we consider in the present paper.
As an unexpected by-product, we show boundedness of strongly spirallike
functions.
\end{abstract}
\thanks{
The first author was supported by Yeungnam University (2010).
The second author was supported in part by JSPS Grant-in-Aid for 
Scientific Research (B) 22340025.
}
\maketitle

\section{Introduction}

Let $\A$ denote the set of analytic functions on the unit disk
$\D=\{z: |z|<1\}$ of the complex plane $\C.$
Set furtheremore $\A_0=\{f\in\A: f(0)=1\}$ and 
$\A_1=\{f\in\A: f(0)=0, f'(0)=1\}.$
We note that a function $h(z)$ belongs to $\A_0$ if and only if
the function $zh(z)$ belongs to $\A_1.$
In what follows, $f(z)/z$ will be regarded as a function in $\A_0$
for $f\in\A_1.$
More concretely, for a function $f(z)=z+a_2z^2+a_3z^3+\cdots$ in $\A_1,$
the function $f(z)/z$ is regarded as the analytic function
$1+a_2z+a_3z^2+\cdots.$
Let $\A_0^\times$ be the set of invertible elements of $\A_0$
with respect to the ordinary multiplication; that is,
$\A_0^\times=\{h\in\A_0: h(z)\ne0,~z\in\D\}.$
In what follows, $\Log h$
means the (single-valued) analytic branch of $\log h$ in $\D$ 
determined by $\Log h(0)=0$ for $h\in\A_0^\times.$
We also set $\Arg h=\Im \Log h$ for $h\in \A_0^\times.$
We note that $\Log$ maps $\A_0^\times$ bijectively onto the complex
vector space $\V=\{f\in\A: f(0)=0\}.$

The set $\es$ consisting of all the univalent funtions in $\A_1$
has been the central object to study in the theory of univalent functions
since early 20th century.

We are interested in classical subclasses of $\es$ in the present paper.
Let us now introduce them.
A function $f\in\A_1$ is called {\it convex} if $f$ maps $\D$
univalently onto a convex domain in $\C.$
We denote by $\K$ the class of convex functions.
It is well known that $f\in\A_1$ is convex if and only if
$$
\Re\left(1+\frac{zf''(z)}{f'(z)}\right)>0,\quad z\in\D.
$$
Let $\lambda$ be a number with $-\pi/2<\lambda<\pi/2.$
For a point $a\ne0$ in $\C,$ the $\lambda$-spiral segment
$[0,a]_\lambda$ is defined to be the set
$\{0\}\cup\{a\exp(-te^{i\lambda}): 0\le t<+\infty\}.$
A domain $\Omega$ in $\C$ is called $\lambda$-spirallike
(about the origin)
if $[0,a]_\lambda\subset\Omega$ for every $a\in\Omega.$
In particular, a $0$-spirallike domain is also called starlike as usual.
A function $f\in\A_1$ is called {\it $\lambda$-spirallike}
if $f$ maps $\D$ univalently onto a $\lambda$-spirallike domain.
The class of $\lambda$-spirallike functions will be denoted by
$\spl(\lambda).$
Set $\spl=\bigcup_{-\pi/2<\lambda<\pi/2}\spl(\lambda).$
The class of starlike functions $\spl(0)$ is also denoted by $\es^*.$
It is also known that $f\in\A_1$ is $\lambda$-spirallike if and only if
$$
\Re\left(e^{-i\lambda}\frac{zf'(z)}{f(z)}\right)>0,\quad
0<|z|<1.
$$
For a real number $\alpha\le1,$ a function $f\in\A_1$ is called
{\it starlike of order $\alpha$} if $\Re(zf'(z)/f(z))\ge\alpha,~z\in\D.$
Let $\es^*(\alpha)$ denote the set of starlike functions of order $\alpha.$
Similarly, for $0<\alpha<1,$ a function $f\in\A_1$ is called 
{\it strongly starlike of order $\alpha$}
if $|\Arg(zf'(z)/f(z))|<\pi\alpha/2,~z\in\D,$
and the set of those functions will be denoted by $\st(\alpha).$

We can extend strong starlikeness to strong spirallikeness in an obvious way.
Let $\lambda\in(-\pi/2,\pi/2)$ and $0<\alpha<1.$
A function $f\in\A_1$ is called {\it strongly $\lambda$-spirallike
of order $\alpha$} if
$$
\left|\Arg\frac{zf'(z)}{f(z)}-\lambda\right|<\frac{\pi\alpha}{2},
\quad z\in\D.
$$
We denote by $\spl(\lambda,\alpha)$ the set of these functions.
When we do not specify $\lambda$ and $\alpha,$ we simply call it
{\it strongly spirallike}.
This sort of classes were first introduced by Bucka and Ciozda \cite{BC73}.

It is an important observation due to Alexander \cite{Alexander15} that
$f(z)$ is convex if and only if $g(z)=zf'(z)$ is starlike.
The mapping $g\mapsto f$ is sometimes called the {\it Alexander transformation}
and will be denoted by $J_1[f]$ in the sequel.
More explicitly,
$$
J_1[f](z)=\int_0^z\frac{f(\zeta)}{\zeta}d\zeta
=\int_0^1f(tz)\frac{dt}{t}
$$
for $f\in\A_1.$
Note also that $J_1(\A_1)=\A_1.$

A function $f\in\A_1$ is called {\it close-to-convex} if
$\Re(e^{-i\lambda}f'/g')>0$ in $\D$ for some 
$g\in\K$ and $\lambda\in(-\pi/2,\pi/2).$
The set of close-to-convex functions will be denoted by $\CC.$

We have the inclusion relations $\K\subset\es^*\subset\CC\subset\es$
and $\es^*\subset\spl\subset\es.$
See \cite{Duren:univ} for basic information about these subclasses
of $\es.$

Several integral operators have been considered by many authors
in connection with univalent functions.
For instance, for $c\in\C,$ we define
$$
I_c[f](z)=\int_0^z\left(f'(\zeta)\right)^c d\zeta
$$
for $f\in\LU=\{f\in\A_1: f'\in\A_0^\times\}$ (`locally univalent'), and
$$
J_c[f](z)=\int_0^z\left(\frac{f(\zeta)}{\zeta}\right)^c d\zeta
$$
for $f\in\ZF=\{f\in\A_1: f(z)/z\in \A_0^\times\}$ (`zero-free' except for
the origin).
Here and hereafter, the complex power $h^c$ for $h\in\A_0^\times$
will be understood as $h^c=\exp(c\,\Log h).$
In particular, we see that $h^c\in\A_0^\times$ for $h\in\A_0^\times$
and $c\in\C.$

Note that $I_c(\LU)\subset\LU$ and $J_c(\ZF)\subset\LU.$
For later convenience, we also set
$\dom_I=\LU, \rg_I=\LU, \dom_J=\ZF$ and $\rg_J=\LU.$

In order to deal with these operators at once, let $X$ represent
one of $I, J$ and $K$ which will be introduced below.
For instance, $X_c$ and $\dom_X$ mean $I_c$ and $\dom_I=\LU,$ respectively,
when $X=I.$

It is an interesting problem to describe or estimate the set
$$
[\M,\N]_X=\{c\in\C: X_c[f]\in\N ~\text{for all}~f\in\M\}
=\{c: X_c(\M)\subset\N\}
$$
for $\M\subset\dom_X$ and $\N\subset\rg_X$
and a family of operators $X_c:\dom_X\to\rg_X, c\in\C.$
(This kind of sets appeared earlier in the authors' paper \cite{KS07AT}.)
When $\M$ consists of a single function $f,$ then we write
$[f,\N]_X$ for $[\{f\},\N]_X.$
We denote by $\D(a,r)=\{z\in\C: |z-a|<r\}$ and by $\bD(a,r)$ its closure.
We summarize known relations of this kind.

\begin{Thm}
\begin{enumerate}
\item
$\bD(0,\tfrac14)\cup\{1\}\subset [\es,\es]_I\subset \bD(0,\tfrac13)\cup\{1\}$
\quad{\rm (Pfaltzgraff \cite{Pf75} and Royster \cite{Roy65})}.
\item
$\bD(0,\tfrac14)\subset [\es,\es]_J\subset \bD(0,\tfrac12)$
\quad{\rm (Y.~J.~Kim and Merkes \cite{KM72})}.
\item
$[\K,\es]_I=[\es^*,\es]_J=\bD(0,\tfrac12)\cup[\tfrac12,\tfrac32]$
\quad{\rm (Aksent'ev and Nezhmetdinov \cite{AN82}, cf.~\cite{KPS02})}.
\item
$[\spl,\es]_J=\bD(0,\tfrac12)$ 
\quad{\rm (Merkes \cite[Corollary 2]{Merkes85})}.
\item
$[\spl(\lambda),\es]_J=\bD(0,\tfrac1{2\cos\lambda})\cup
[\tfrac{e^{-i\lambda}}{2\cos\lambda},\tfrac{3e^{-i\lambda}}{2\cos\lambda}]$
\quad{\rm (\cite{KS07AT})}.
\end{enumerate}
\end{Thm}

In the present paper, we would like to propose yet another operator
$K_c$ for $c\in\C$ defined by
$$
K_c[f](z)=z\left(\frac{f(z)}{z}\right)^c
$$
for $f\in\ZF.$
This will be called the {\it power deformation} of $f$ with exponent $c.$
Let $\dom_K=\rg_K=\ZF.$
Of course, the present paper is not the first to define it.
Indeed, this simple operation was used at many places before
(for instance, \cite{Schild58}, \cite{Pinchuk71}, \cite{Merkes85}).
It seems, however, that the operators $K_c$ have not been studied
systematically in the literature.

Introduction of this operator is motivated by the following facts:
\begin{equation}\label{eq:spl}
K_{e^{i\lambda}\cos\lambda}(\es^*)=\spl(\lambda),
\quad -\frac{\pi}2<\lambda<\frac{\pi}2
\end{equation}
(see \cite{Merkes85}) and
\begin{equation}\label{eq:order}
K_{1-\alpha}(\es^*)=\es^*(\alpha),\quad 0\le \alpha<1
\end{equation}
(see \cite{Pinchuk71}, \cite{Merkes85}).
These relations easily follow from the relation
\begin{equation}\label{eq:fc}
\frac{zf_c'(z)}{f_c(z)}=1-c+c\frac{zf'(z)}{f(z)},
\end{equation}
where $f_c=K_c[f].$

Thus, several typical subclasses of $\es$ can be obtained as
power deformations of $\es^*.$

We will show the following relations.

\begin{thm}\label{thm:main}
\hspace{1cm}
\begin{enumerate}
\item\label{item:st2s}
$[\es^*,\es]_K=[\es^*,\spl]_K=\bD(\tfrac12,\tfrac12).$
\medskip
\item\label{item:sta2s}
$[\es^*(\alpha),\es]_K=[\es^*(\alpha),\spl]_K
=\bD\big(\frac1{2(1-\alpha)},\frac1{2(1-\alpha)}\big)$
for $0\le\alpha<1.$
\medskip
\item\label{item:k2s}
$[\K,\es]_K=[\K,\spl]_K=\bD(1,1).$
\medskip
\item\label{item:spl2s}
$[\spl(\lambda),\es]_K=[\spl(\lambda),\spl]_K
=\bD\big(\frac{1-i\tan\lambda}2, \frac1{2\cos\lambda}\big)$
for $\lambda\in(-\pi/2,\pi/2).$
\medskip
\item\label{item:sp2s}
$[\spl,\es]_K=[\spl,\spl]_K=[0,1].$
\medskip
\item\label{item:sst2s}
$[\st(\alpha),\es]_K=[\st(\alpha),\spl]_K=
\bD\big(\frac{1-i\cot\frac{\pi\alpha}2}2, \frac1{2\sin\frac{\pi\alpha}2}\big)
\cup
\bD\big(\frac{1+i\cot\frac{\pi\alpha}2}2, \frac1{2\sin\frac{\pi\alpha}2}\big)$
\quad for $0<\alpha<1.$
\medskip
\item\label{item:ssp2s}
$[\spl(\lambda,\alpha),\es]_K=[\spl(\lambda,\alpha),\spl]_K=
\bD\Big(\frac{1-i\tan\lambda_+}2, \frac1{2\cos\lambda_+}\Big)\cup
\bD\Big(\frac{1-i\tan\lambda_-}2, \frac1{2\cos\lambda_-}\Big)$
\quad for $|\lambda|<\pi\alpha/2<\pi/2,$ where
$\lambda_\pm=\lambda\pm\pi(1-\alpha)/2.$
\medskip
\item\label{item:s2s}
$[\es,\es]_K
=[\CC,\es]_K=\{0,1\}.$
\end{enumerate}
\end{thm}


As an application of our investigation of power deformations,
we obtain the following result, which is used in the second author's
paper \cite{Sugawa10SSP}.

\begin{thm}\label{thm:bdd}
Let $f$ be a strongly spirallike function.
Then $\Log f(z)/z$ is bounded on $\D.$
In particular, $f(z)$ is bounded on $\D.$
\end{thm}

We note that boundedness of strongly starlike functions 
is due to Brannan and Kirwan \cite{BK69}.

\section{Fundamental facts}

In this section, we collect fundamental properties of the operators
$I_c, J_c, K_c$ and the sets $[\M,\N]_X$ of exponents for $X=I,J,K.$

We first observe that the Alexander transformation $J_1$ maps the class $\ZF$
of zero-free functions onto 
$\LU,$ the class of locally univalent functions, in a one-to-one manner.
By definition, we have
$$
J_c=I_c\circ J_1=J_1\circ K_c
$$
for $c\in\C.$
In particular, we have $K_c=J_1\inv\circ I_c\circ J_1.$
Furthermore, Alexander's observation gives $J_1(\es^*)=\K.$
Therefore, we have $J_c(\es^*)=I_c(\K)$ for $c\in\C.$

Recall now that the set $\V=\{f\in\A: f(0)=0\}$ is a subspace
of the complex vector space $\A.$
We consider the bijective maps $\Phi:\LU\to \V$ and 
$\Psi:\ZF\to \V$ defined by $\Phi[f]=\Log f'$ and $\Psi[f](z)=\Log f(z)/z.$
Then the operators $I_c$ and $K_c$ can be viewed as
scalar multiplication in $\V$ when we identify $\LU$ and $\ZF$ with $\V$ through
the maps $\Phi$ and $\Psi,$ respectively.
In other words, $I_c[f]=\Phi\inv(c\Phi[f])$
and $K_c[f]=\Psi\inv(c\Psi[f]).$
In particular, we easily have the relations
$I_c\circ I_{c'}=I_{cc'}$ and
$K_c\circ K_{c'}=K_{cc'}$ for $c,c'\in\C.$
Moreover, we can even introduce linear structures to the sets
$\LU$ and $\ZF,$ although we will not go into details in the present paper.
Indeed, such a linear structure on $\LU$ was first considered by
Hornich \cite{Hor69} (see also \cite{KPS02}).

We now collect obvious properties of the sets $[\M,\N]_X$
for $X=I,J,K.$

\begin{lem}\label{lem:basic}
Let $X$ represent one of $I, J, K$ and
let $\M, \M', \M_\lambda\subset\dom_X~(\lambda\in\Lambda), \N, \N'\subset\rg_X.$
Then the following hold:
\begin{enumerate}
\item\label{item:b1}
$[\M,\N]_X\supset [\M',\N]_X$ if $\M\subset \M'.$
\medskip
\item\label{item:b2}
$[\M,\N]_X\subset [\M,\N']_X$ if $\N\subset \N'.$
\medskip
\item\label{item:b3}
$[\bigcup_{\lambda\in\Lambda}\M_\lambda,\N]_X
=\bigcap_{\lambda\in\Lambda}[\M_\lambda,\N]_X.$
\medskip
\item\label{item:b4}
$[\bigcap_{\lambda\in\Lambda}\M_\lambda,\N]_X
\supset\bigcup_{\lambda\in\Lambda}[\M_\lambda,\N]_X.$
\medskip
\item\label{item:c}
$[X_c(\M),\N]_X=\frac1c\,[\M,\N]_X$ for $c\in\C\setminus\{0\}$ and for $X=I,K.$
\medskip
\item\label{item:b6}
$[\M,\N]_X$ is a closed subset of $\C$ if
$\N$ is closed in the topology of local uniform convergence on $\D.$
\medskip
\item\label{item:b7}
$[\M,\N]_K=[J_1(\M),J_1(\N)]_I.$
\end{enumerate}
\end{lem}

Here, we define $c E=\{cz: z\in E\}$ for $E\subset\C$ and $c\in\C.$
We remark that $\es, \es^*(\alpha), \K, \CC,$
$\SS(\alpha), \spl(\lambda), \spl(\lambda,\alpha), \spl, \LU, \ZF$ 
are all closed in the topology of local uniform convergence on $\D.$

The power deformation effects on boundedness.
We summarize a few facts about it.

\begin{lem}\label{lem:unbdd}
For a function $f\in \ZF$ and $c\in\C,$ let $f_c=K_c[f].$
\begin{enumerate}
\item\label{item:unbdd1}
If $\Log f(z)/z$ is bounded in $\D,$ then so is $\Log f_c(z)/z$ for
every $c\in\C.$
\item
If $\log|f(z)/z|$ is unbounded in $\D,$ then so is
$\log|f_c(z)/z|$ for every $c>0.$
\item\label{item:univ}
Suppose that $f$ is unbounded and univalent in $\D$
and that $\Arg f(z)/z$ is bounded in $\D.$
Then $f_c$ is never univalent when $\Re c<0$ while
$f_c$ is unbounded when $\Re c>0.$
\end{enumerate}
\end{lem}

\begin{pf}
Assertions (1) and (2) are clear when we look at the relation
$\Log f_c(z)/z=c\,\Log f(z)/z.$

We prove assertion (3).
Let $c=a+ib.$
By assumption, we have a sequence $z_n~(n=1,2,\dots)$ in $\D$
such that $|f(z_n)|\to\infty$ and $|z_n|\to1$ as $n\to\infty.$
Since we have the relation
$$
\log\left|\frac{f_c(z)}{z}\right|=a\log\left|\frac{f(z)}{z}\right|
-b \,\Arg\frac{f(z)}{z},
$$
$|f_c(z_n)/z_n|\to 0$ as $n\to\infty$ if $b<0.$
Then, $f_c$ is never univalent.
Also, the above relation tells us that $f_c$ is unbounded if $b>0.$
\end{pf}

For a subclass $\M$ of $\ZF,$ we denote by $V(\M)$
the variability region of the quantity $zf'(z)/f(z)$ for $f\in\M;$
more concretely,
$$
V(\M)=\{zf'(z)/f(z): f\in\M, z\in\D\}.
$$
Note that $V(\M)$ is a domain (a connected non-empty open set)
unless $\M\subset\{\id\}.$
This has a close connection with $[\M,\es]_K.$
Let $T$ be the M\"obius transformation defined by 
$$
T(w)=\frac1{1-w}.
$$
Then, we have $[\M,\es]_K\subset \C\setminus T(V(\M))$ by
the following result.

\begin{lem}\label{lem:T}
For a subclass $\M$ of $\ZF,$ the set $[\M,\LU]_K$
and the variability region $V(\M)$ of $zf'(z)/f(z)$ are related by
$$
[\M,\LU]_K=\C\setminus T(V(\M)).
$$
\end{lem}

\begin{pf}
Let $c$ be a finite point in $T(V(\M)).$
Then there are $f\in\M$ and $z_0\in\D$ such that $c=T(z_0f'(z_0)/f(z_0));$
namely, $z_0f'(z_0)/f(z_0)=1-1/c.$
Then by \eqref{eq:fc} the function $f_c=K_c[f]$ satisfies 
$$
\frac{z_0f_c'(z_0)}{f_c(z_0)}=1-c+c\frac{z_0f'(z_0)}{f(z_0)}
=0.
$$
In particular, $K_c[f]$ is not locally univalent at $z_0$
and therefore $c\notin[\M,\LU]_K.$

We can also trace back the above argument to prove the converse.
\end{pf}

For an $f\in\ZF,$ set $V(f)=\{zf'(z)/f(z): z\in\D\}.$
Then, in particular, we have the relation
$$
[f,\LU]_K=\C\setminus T(V(f)).
$$

We can also derive the following corollary.

\begin{cor}\label{cor:bdd}
Let $\M$ be a subclass of $\ZF$ which contains a function $f\ne\id.$
Then $[\M,\LU]_K$ is a compact subset of $\C.$
\end{cor}

\begin{pf}
Since $V(\M)$ is a domain containing $1,$ the image
$T(V(\M))$ under $T$ is a domain in the Riemann sphere containing $\infty.$
Therefore, its complement $[\M,\LU]_K$ is compact in $\C.$
\end{pf}

Therefore, $[\M,\N]_K$ is compact when $\M$ and $\N$ are chosen from
$\es, \es^*(\alpha), \K, \CC, \spl(\lambda,\alpha), \spl.$

We summarize information about the variabirity regions of typical
subclasses of $\es.$

\begin{lem}\label{lem:vr}
One has the following relations:
\begin{enumerate}
\item
$V(\es^*)=\{w: \Re w>0\}.$
\item
$V(\es^*(\alpha))=\{w: \Re w>\alpha\}.$
\item
$V(\K)=\{w: \Re w>1/2\}.$
\item
$V(\spl(\lambda))=\{w: \Re e^{-i\lambda}w>0\}.$
\item
$V(\spl)=\C\setminus(-\infty,0].$
\item
$V(\st(\alpha))=\{w: |\arg w|<\pi\alpha/2\}.$
\item
$V(\spl(\lambda, \alpha))=\{w: |\arg w-\lambda|<\pi\alpha/2\}.$
\item
$V(\es)=V(\CC)=\C\setminus\{0\}.$
\end{enumerate}
\end{lem}

\begin{pf}
We have to show a relation of the form $V(\M)=B$ for a class $\M$
and a subdomain $B$ of $\C$ in each case.
When $V(\M)\subset B$ is trivial by the definition of $\M,$
we just give a function $f\in\M$ such that
$zf'(z)/f(z)$ covers the domain $B$ in order to show $B\subset V(\M).$

\noindent
(1) Consider the Koebe function $k(z)=z/(1-z)^2.$ 

\noindent
(2) Consider the function $K_{1-\alpha}[k](z)=z/(1-z)^{2(1-\alpha)}.$

\noindent
(3) E.~Strohh\"acker showed the relation $\K\subset\es^*(\frac12)$
(see \cite[p.~251]{Duren:univ} for instance).
Therefore, we have $V(\K)\subset\{w: \Re w>1/2\}.$
On the other hand, $l(z)=z/(1-z)$ is convex and $zl'(z)/l(z)=1/(1-z)$
maps $\D$ conformally onto the half-plane $\Re w>1/2.$
Therefore, we have $V(\K)=\{w: \Re w>1/2\}.$

\noindent
(4) Consider the function $K_{e^{i\lambda}\cos\lambda}[k](z)
=z/(1-z)^{2e^{i\lambda}\cos\lambda}.$

\noindent
(5) This is clear because $V(\spl)=\cup_{\lambda} V(\spl(\lambda)).$

\noindent
(6) Consider the function $f\in\A_1$ determined by $zf'(z)/f(z)
=(\frac{1+z}{1-z})^\alpha.$

\noindent
(7) Consider the function $f\in\A_1$ determined by $zf'(z)/f(z)
=(\frac{1+ze^{2i\lambda/\alpha}}{1-z})^\alpha.$

\noindent
(8) The assertion $V(\CC)=\C\setminus\{0\}$ can be found in \cite{Wang11}.
Since $V(\CC)\subset V(\es)\subset\C\setminus\{0,1\},$ 
the other assertion follows, too.
\end{pf}

\section{Proof of main results}

\begin{pf}[Proof of Theorem \ref{thm:main}]

We need to prove the assertion $[\M,\es]_K=[\M,\spl]_K=A$
for $\M=\es^*, \es^*(\alpha), \K, \spl(\lambda),$
$\spl, \st(\alpha),\spl(\lambda,\alpha), \es, \CC$ 
and the subset $A\subset\C$ which
appears in the right-hand side of the relation in the corresponding assertion
(though we should omit $[\M,\spl]_K$ in the case of \eqref{item:s2s}).
First we observe that the set $A$ is indeed equal to $\C\setminus T(V(\M))$
in each case by virtue of Lemma \ref{lem:vr}.
Therefore, by Lemma \ref{lem:T} and Lemma \ref{lem:basic} \eqref{item:b2},
we obtain
$$
[\M,\spl]_K\subset [\M,\es]_K\subset[\M,\LU]_K=\C\setminus T(V(\M))=A.
$$
Therefore, it is enough to show that $A\subset [\M,\spl]_K$
with the exception of \eqref{item:s2s}.
We will take this strategy unless a simpler way is available.
We divide the proof into several pieces according to the numbering.

\medskip\noindent
[Proof of \eqref{item:st2s}:]
We show the implication $\bD(\frac12,\frac12)\subset[\es^*,\spl]_K.$
Let $f\in\es^*$ and set $f_c=K_c[f]$ for $c\in\C.$
Then, by \eqref{eq:order}, we have $f_c\in\es^*(1-c)\subset\es^*$
for $0\le c\le 1.$
Next, by \eqref{eq:spl}, we see that
$$
f_{ce^{i\lambda}\cos\lambda}=K_{e^{i\lambda}\cos\lambda}[f_c]
\in \spl(\lambda)\subset\spl
$$
for $\lambda\in(-\pi/2,\pi/2).$
In view of the relation $e^{i\lambda}\cos\lambda=(e^{2i\lambda}+1)/2,$
we obtain 
$$
\{ce^{i\lambda}\cos\lambda: 0\le c\le 1, -\pi/2<\lambda<\pi/2\}
=\bD\big(\tfrac12,\tfrac12\big).
$$
Thus we have shown that $\bD(\frac12,\frac12)\subset[\es^*,\spl]_K.$

\medskip\noindent
[Proof of \eqref{item:sta2s} and \eqref{item:spl2s}:]
We combine Lemma \ref{lem:basic} \eqref{item:c} with \eqref{eq:order}
and \eqref{eq:spl} to obtain \eqref{item:sta2s} and \eqref{item:spl2s}.
Here, we note the relation $1/(e^{i\lambda}\cos\lambda)=1-i\tan\lambda.$

\medskip\noindent
[Proof of \eqref{item:k2s}:]
By the Strohh\"acker theorem: $\K\subset\es^*(\frac12)$ which is mentioned
in the proof of Lemma \ref{lem:vr},
we obtain 
$$
[\K,\spl]_K\supset [\es^*(\tfrac12),\spl]_K=\bD(1,1).
$$

\medskip\noindent
[Proof of \eqref{item:sp2s}:]
It is enough to show that $[0,1]\subset[\spl,\spl]_K.$
This follows from the fact that $[0,1]\subset[\spl(\lambda),\spl]_K$
for every $\lambda\in(-\pi/2,\pi/2)$ by \eqref{item:spl2s}.

\medskip\noindent
[Proof of \eqref{item:sst2s}:]
Since $\st(\alpha)=\spl(0,\alpha),$ this follors from \eqref{item:ssp2s}.

\medskip\noindent
[Proof of \eqref{item:ssp2s}:]
Since $\spl(\lambda,\alpha)=\spl(\lambda_+)\cap\spl(\lambda_-),$
Lemma \ref{lem:basic} \eqref{item:b4} yields the relation
$$
[\spl(\lambda,\alpha),\spl]_K\supset
\bD\Big(\frac{1-i\tan\lambda_+}2, \frac1{2\cos\lambda_+}\Big)\cup
\bD\Big(\frac{1-i\tan\lambda_-}2, \frac1{2\cos\lambda_-}\Big).
$$

\medskip\noindent
[Proof of \eqref{item:s2s}:]
It is enough to see $\{0,1\}\subset [\es,\es]_K.$
This is trivial.
\end{pf}

\begin{rem}
As we saw in the proof, we actually showed the relations
$$
[\M,\es]_K=[\M,\LU]_K=\C\setminus T(V(\M))
$$
for $\M=\es^*, \es^*(\alpha), \K, \spl(\lambda), \spl, \st(\alpha), 
\spl(\lambda,\alpha), \CC, \es.$
Under this situation, if a function $f_0\in\M$ satisfies $V(f_0)=V(\M),$
then $\C\setminus T(V(\M))=[\M,\es]_K\subset [f_0,\LU]_K=\C\setminus T(V(f_0))$
and therefore $[f_0,\es]_K=[\M,\es]_K.$
The class $\es$ can be replaced by any class as long as the above relations are
valid.
For instance, the Koebe function $k$ satisfies
$[k,\es]_K=[k,\spl]_K=[k,\LU]_K=\bD(\frac12,\frac12).$
\end{rem}

In order to prove Theorem \ref{thm:bdd}, we recall the following result.

\begin{lem}[Goodman \cite{Good53}]\label{lem:good}
$|\Arg f(z)/z|\le 2\arcsin|z|<\pi,~|z|<1,$ for $f\in\es^*.$
\end{lem}

We are now ready to prove Theorem \ref{thm:bdd}.

\begin{pf}[Proof of Theorem \ref{thm:bdd}]
Let $f$ be strongly $\lambda$-spirallike of order $\alpha$
with $|\lambda|<\pi\alpha/2<\pi/2.$
Put $g=K_{e^{-i\lambda}/\cos\lambda}[f].$
By Theorem \ref{thm:main} \eqref{item:ssp2s} together
with Lemma \ref{lem:basic} \eqref{item:c}, we have
\begin{align*}
[g,\es]_K
&=[K_{e^{-i\lambda}/\cos\lambda}[f],\es]_K
=e^{i\lambda}\cos\lambda[f,\es]_K \\
&\supset e^{i\lambda}\cos\lambda\left(
\bD\Big(\frac{1-i\tan\lambda_+}2, \frac1{2\cos\lambda_+}\Big)\cup
\bD\Big(\frac{1-i\tan\lambda_-}2, \frac1{2\cos\lambda_-}\Big)\right),
\end{align*}
where $\lambda_\pm=\lambda\pm\pi(1-\alpha)/2.$
Observe that $[g,\es]_K$ is not contained in the closed right half-plane.

On the other hand, by \eqref{eq:spl}, $g\in\es^*$ because $f\in\spl(\lambda).$
Thus $g$ is univalent and $\Arg g(z)/z$ is bounded by Lemma \ref{lem:good}.
We now suppose that $g$ was unbounded in $\D.$
Then Lemma \ref{lem:unbdd} impliles that $[g,\es]_K$ would be contained in
the closed right half-plane $\Re c\ge 0.$
This is a contradiction.
We have shown that $g$ is bounded, and hence,
$\Log g(z)/z$ is bounded.
We now have boundedness of $\Log f(z)/z$
by Lemma \ref{lem:unbdd} \eqref{item:unbdd1}.
\end{pf}

It is somewhat strange that we obtained a boundedness result for strongly
spirallike functions without making any
concrete estimate of functions involved.
We also note that the above $g$ satisfies the relation
$zg'(z)/g(z)=c zf'(z)/f(z)+1-c,$ where $c=e^{-i\lambda}/\cos\lambda
=1-i\tan\lambda.$
Therefore, $g$ is not necessarily strongly starlike unless $\lambda=0.$

\def\cprime{$'$} \def\cprime{$'$} \def\cprime{$'$}
\providecommand{\bysame}{\leavevmode\hbox to3em{\hrulefill}\thinspace}
\providecommand{\MR}{\relax\ifhmode\unskip\space\fi MR }
\providecommand{\MRhref}[2]{%
  \href{http://www.ams.org/mathscinet-getitem?mr=#1}{#2}
}
\providecommand{\href}[2]{#2}

\end{document}